\documentclass[a4paper]{article}
\usepackage[numbers]{natbib}
\usepackage{amsthm}
\usepackage{amsmath,amssymb}
\usepackage[scale=.7]{geometry}

\theoremstyle{plain}
\newcounter{theoremCounter}
\numberwithin{theoremCounter}{section}

\theoremstyle{definition}

\newtheorem{remark}[theoremCounter]{Remark}

\usepackage[usenames,dvipsnames]{xcolor}
\usepackage{graphicx}        
\usepackage[bottom]{footmisc}

\usepackage{pst-all}
\usepackage{pstricks-add}
\usepackage{psfrag}

\usepackage{hyperref}

\makeindex             

\newcommand{\fett}[1]{\mbox{\boldmath$#1$}} 
\newcommand{\tena}[1]{\lowercase{\fett{#1}}} 
\newcommand{\tenb}[1]{\uppercase{\fett{#1}}} 
\newcommand{\dyad}{\otimes} 

\begin{document}
\makeatletter\let\@fnsymbol\@arabic\makeatother 
\title{On the dislocation density tensor in the Cosserat theory of elastic shells}
\author{Mircea B\^{\i}rsan\thanks{%
Lehrstuhl f\"{u}r Nichtlineare Analysis und Modellierung, Fakult\"{a}t f\"{u}r Mathematik,
Universit\"{a}t Duisburg-Essen, Thea-Leymann Str. 9, 45127 Essen, Germany; and  Alexandru Ioan Cuza University of Ia\c si, Department of Mathematics,  Blvd.
Carol I, no. 11, 700506 Ia\c si, Romania, \url{mircea.birsan@uni-due.de}}%
\quad and\quad%
Patrizio Neff\thanks{%
Head of Lehrstuhl f\"{u}r Nichtlineare Analysis und Modellierung, Fakult\"{a}t f\"{u}r
Mathematik, Universit\"{a}t Duisburg-Essen,  Thea-Leymann Str. 9, 45127 Essen, Germany \url{patrizio.neff@uni-due.de}}%
}
%
%
\maketitle


\abstract{We consider the Cosserat continuum in its finite strain setting and discuss the dislocation density tensor as a possible alternative curvature strain measure in three-dimensional Cosserat models and in Cosserat shell models. We establish a close relationship (one-to-one correspondence) between the new shell dislocation density tensor and the bending-curvature tensor of 6-parameter shells.}

\section{Indroduction}
\label{sec:1}

The Cosserat type theories have recently seen a tremendous renewed interest for their prospective applicability to model physical effects beyond the classical ones. These comprise notably the so-called size--effects (``smaller is stiffer'').

In a finite strain Cosserat-type framework the group of proper rotations $\mathrm{SO}(3)$ has a dominant place. The original idea of the Cosserat brothers \cite{Cosserat09} to consider independent rotational degrees of freedom in addition to the macroscopic displacement was heavily motivated by their treatment of plate and shell theory. Indeed, in shell theory it is natural to attach a preferred orthogonal frame (triad) at any point of the surface, one vector of which is the normal to the midsurface, the other two vectors lying in the tangent plane.
This is the notion of the ``tri\`edre cach\'e''. The idea to consider then an orthogonal frame which is not strictly linked to the surface, but constitutively coupled, leads to the notion of the ``tri\`edre mobile''. And this then is already giving rise to a prototype Cosserat shell (6-parameter) theory. For an insightful review of various Cosserat-type shell models we refer to \cite{Altenbach-Erem-Review}.

However, the Cosserat brothers have never proposed any more specific constitutive framework, apart from postulating euclidean invariance (frame-indifference) and hyperelasticity. For specific problems it is necessary to choose a constitutive framework and to determine certain strain and curvature measures. This task is still not conclusively done, see e.g.\\\ \cite{Pietraszkiewicz09}.

Among the existing models for Cosserat-type shells we mention the theory of simple elastic shells \cite{Altenbach04}, which has been developed by \cite{Zhilin76,Zhilin06} and \cite{Altenbach-Zhilin-82,Altenbach-Zhilin-88}. Later, this theory has been successfully applied to describe the mechanical behavior of laminated, functionally graded, viscoelastic or porous plates in \cite{Altenbach00,Altenbach-Erem08,Altenbach-Erem-AM-09,Altenbach-Erem-10} and of multi-layered, orthotropic, thermoelastic shells in \cite{Birsan-Alten-MMAS-10,Birsan-Alten-2011,Birsan-Sadow-JTS-13,Birsan-ACME-15}.
Another remarkable approach is the general 6-parameter theory of elastic shells presented in \cite{Libai98,Pietraszkiewicz-book04,Pietraszkiewicz04}. Although the starting point is different, one can see that the kinematical structure of the nonlinear 6-parameter shell theory is identical to that of a Cosserat shell model, see also \cite{Birsan-Neff-MMS-2014,Birsan-Neff-L54-2014}.

In this paper we would like to draw attention to alternative curvature measures, motivated by dislocation theory, which can also profitably be used in the three-dimensional Cosserat model and the Cosserat shell model. The object of interest is Nye's dislocation density tensor $\,\mathrm{Curl}\,\tenb P\,$. Within the restriction to proper rotations it turns out that Nye's tensor provides a complete control of all spatial derivatives of rotations \cite{Neff_curl08} and we rederive this property for micropolar continua using general curvilinear coordinates.
Then we focus on shell-curvature measures and define a new shell dislocation density tensor using the surface Curl operator. Then, we prove that a relation analogous to Nye's formula holds also for Cosserat (6-parameter) shells.

The paper is structured as follows. In Section \ref{sec:2} we present the kinematics of a three-dimensional Cosserat continuum, as well as the appropriate strain measures and curvature strain measures, written in curvilinear coordinates.
Here, we show the close relationship between the  wryness tensor and the dislocation density tensor, including the corresponding Nye's formula. In Section \ref{sec:3} we define the $\,\mathrm{Curl}\,$ operator on surfaces and present several representations using surface curvilinear coordinates. These relations are then used in Section \ref{sec:4} to introduce the new shell dislocation density tensor and to investigate its relationship to the elastic shell bending-curvature tensor of 6-parameter shells.

\section{Strain measures of a three-dimensional Cosserat model in curvilinear coordinates}
\label{sec:2}

Let $\cal B$ be a Cosserat elastic body which occupies in its reference (initial) configuration the domain $\Omega_\xi\subset\mathbb{R}^3$.  A generic point of $\Omega_\xi$ will be denoted by $(\xi_1,\xi_2,\xi_3)$. The deformation of the Cosserat body is described by a vectorial map $\tena{\varphi}_\xi$  and  a microrotation tensor $\tenb{R}_\xi\,$,
\begin{align*}
\tena\varphi_\xi:\Omega_\xi \rightarrow\Omega_c\,, \qquad \tenb{R}_\xi:\Omega_\xi \rightarrow \mathrm{SO}(3)\, ,
\end{align*}
where $\Omega_c\,$ is the deformed (current) configuration. Let $(x_1,x_2,x_3)$ be some general curvilinear coordinates system on $\Omega_\xi\,$. Thus, we have a parametric representation  $\, \tena\Theta \,$ of the domain $\,\Omega_\xi$
\begin{align*}
\tena\Theta:\Omega \rightarrow\Omega_\xi\,, \qquad \tena\Theta(x_1,x_2,x_3)=(\xi_1,\xi_2,\xi_3),
\end{align*}
where   $\,\Omega\subset\mathbb{R}^3$ is a bounded domain with Lipschitz boundary $\partial\Omega$.
The covariant base vectors with respect to these curvilinear coordinates are denoted by  $\tena g_i$ and the contravariant base vectors by $\tena g^j$ ($i,j=1,2,3$), i.e.
\begin{align*}
    \tena g_i=\dfrac{\partial\tena\Theta }{\partial x_i}\,=\tena\Theta_{,i}\,,\qquad
    \tena g^j\cdot \tena g_i=\delta^j_i\,,
\end{align*}
where $\delta^j_i$ is the Kronecker symbol. We employ the usual conventions for indices: the Latin indices $i,j,k,...$ range over the set $\{1,2,3\}$, while the Greek indices $\alpha,\beta,\gamma,...$ are confined to the range $\{1,2\}\,$; the comma preceding an index $i$ denotes partial derivatives with respect to $x_i\,$; the Einstein summation convention over repeated indices is also used.

Introducing the deformation function $\,\tena \varphi\,$   by the composition
\begin{align*}
    \tena \varphi:=\tena \varphi_\xi\circ\tena\Theta :\Omega \rightarrow\Omega_c\,,\qquad \tena \varphi(x_1,x_2,x_3): = \tena \varphi_\xi\big( \tena\Theta(x_1,x_2,x_3)\big),
\end{align*}
we can express the (elastic) deformation gradient $\tenb F$ as follows
\begin{align*}
 \tenb F := \nabla_\xi\,\tena\varphi_\xi ( \xi_1,\xi_2,\xi_3)\, =\nabla_x\,\tena \varphi(x_1,x_2,x_3)\cdot\big[\nabla_x \tena \Theta(x_1,x_2,x_3)\big]^{-1}\,.
\end{align*}
Using the direct tensor notation, we can write
\begin{align*}
  \nabla_x\tena\varphi=\tena\varphi_{,i}\dyad \tena e_{i}\,,\qquad \nabla_x\tena\Theta=\tena g_{i}\dyad \tena e_{i}\,,\qquad \big[\nabla_x\,\tena\Theta]^{-1}= \tena e_{j}\dyad \tena g^{j},
\end{align*}
where $\,\tena e_i\,$ are the unit vectors along the coordinate axes $Ox_i$ in the parameter domain $\,\Omega\,$. Then, the deformation gradient can be expressed by
\begin{align*}
  \tenb F = \tena\varphi_{,i}\dyad \tena g^{i}.
\end{align*}

\begin{figure}
\centering
\includegraphics[scale=0.7]{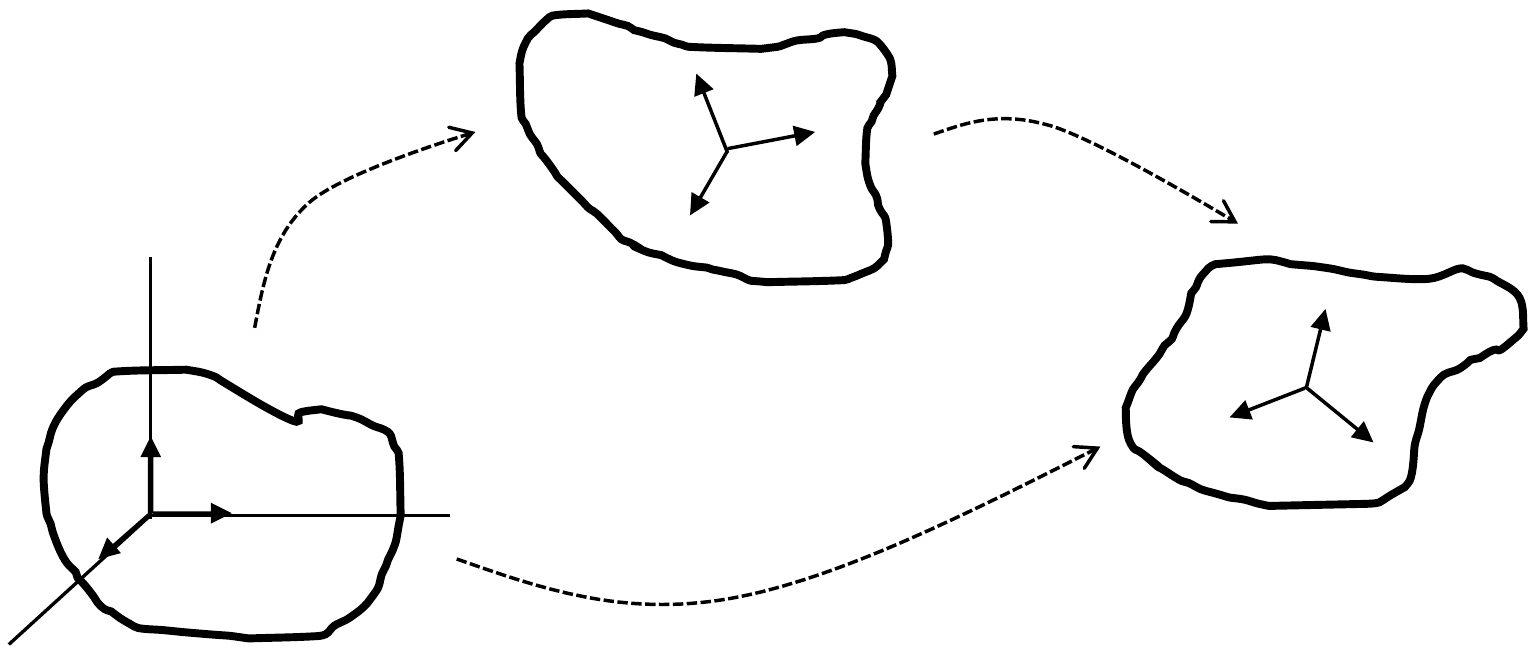}
\put(-255,9){$\Omega$} \put(-312,8){$x_1$} \put(-278,77){$x_3$} \put(-282,20){\small{$O$}} \put(-225,31){$x_2$}
\put(-288,14){$\tena e_1$} \put(-270,32){$\tena e_2$} \put(-291,40){$\tena e_3$}
\put(-169,83){$\tena d^0_1$} \put(-150,95){$\tena d^0_2$} \put(-180,112){$\tena d^0_3$}
\put(-65,39){$\tena d_1$} \put(-39,35){$\tena d_2$} \put(-40,62){$\tena d_3$}
\put(-170,18){$\overline{\tenb R}(x_i)$} \put(-165,3){$ \tena \varphi(x_i)$}
\put(-265,98){$\tenb Q_0(x_i)$} \put(-248,83){$ \tenb \Theta(x_i)$}
\put(-110,111){$\tenb R_\xi(\xi_i)=\tenb Q_e(x_j)$} \put(-108,92){$ \tena \varphi_\xi(\xi_i)$}
\put(-165,65){$\Omega_\xi$}  \put(-48,21){$\Omega_c$} 
\caption{The reference (initial) configuration $\Omega_\xi$ of the Cosserat continuum, the deformed (current) configuration $\Omega_c$ and the parameter domain $\Omega$ of the curvilinear coordinates $(x_1,x_2,x_3)$. The triads of directors $\{ \tena d_i\}$  and  $\{ \tena d^0_i\}$ satisfy the relations $\tena{d}_i=\tenb{Q}_e\tena{d}_i^0=  \overline{\tenb R}\tena{e}_i$  and $ \tena{d}_i^0=  \tenb Q_0\tena{e}_i\,$, where $\tenb{Q}_e$ is the elastic microrotation field, $\tenb{Q}_0$   the initial microrotation, and $\overline{\tenb R}$ the total microrotation field.}
\label{fig:1}       
\end{figure}
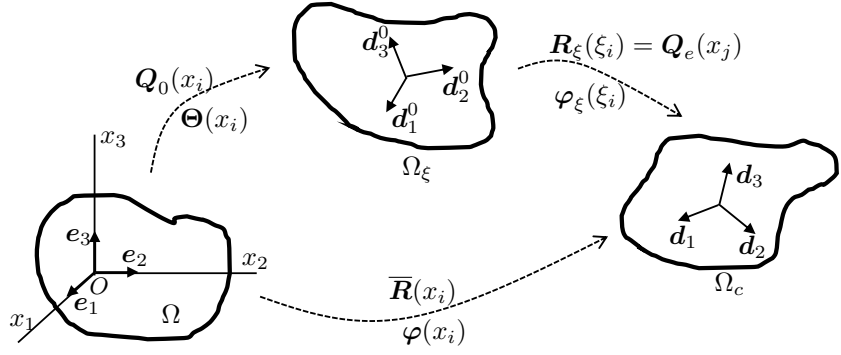

The orientation and rotation of points in Cosserat (micropolar) media can also be described by means of triads of orthonormal vectors (called \emph{directors}) attached to every point. We denote by $\{\tena{d}_i^0 \}$  the triad of directors ($i=1,2,3$) in the reference configuration $\,\Omega_\xi\,$ and by $\{\tena{d}_i \}$ the directors in the deformed configuration $\,\Omega_c\,$, see Figure 1. We introduce the  \emph{elastic microrotation} $\tenb Q_e$ as the composition
\begin{align*}
     \tenb Q_e:= \tenb R_\xi\circ\tena\Theta :\Omega \rightarrow \mathrm{SO}(3 ),\qquad \tenb Q_e (x_1,x_2,x_3):= \tenb R_\xi\big(\tena\Theta(x_1,x_2,x_3)\big),
\end{align*}
which can be characterized with the help of the directors by the relations
$$ \,\tenb{Q}_e\tena{d}_i^0=\tena{d}_i\,, \quad\text{i.e.}, \qquad \tenb Q_e=\tena d_i\dyad\tena d_i^0\,.$$
Let $\tenb Q_0$ be the \emph{initial microrotation} (describing the position of the directors in the reference configuration $\,\Omega_\xi\,$)
$$ \,\tenb{Q}_0\tena{e}_i =\tena{d}_i^0\,, \quad\text{i.e.}, \qquad \tenb Q_0=\tena d_i^0\dyad\tena e_i\,.$$
Then, the \emph{total microrotation} $\overline{\tenb R}$ is given by
\begin{align*}
    \overline{\tenb R}:\Omega \rightarrow \mathrm{SO}(3 ),\qquad \overline{\tenb R}(x_i):=\tenb Q_e(  x_i)\,\tenb Q_0(  x_i) = \tena d_i(  x_i)\dyad\tena e_i\,.
\end{align*}
The  non-symmetric Biot-type stretch tensor (the elastic \emph{first Cosserat deformation tensor}, see \cite{Cosserat09}, p. 123, eq. (43)) is now
\begin{align*}
    \overline{\tenb{U}}_e:=\tenb Q_e^T\, \tenb F= \big(\tena d_i^0\dyad\tena d_i\big)\,\big(\tena \varphi,_j\dyad \tena g^j \big)= \big(\tena \varphi,_j\cdot\, \tena d_i\big)\,
    \tena d_i^0\dyad \tena g^j.
\end{align*}
and the non-symmetric \emph{strain tensor} for nonlinear micropolar materials is defined by
\begin{align*}
    \overline{\tenb{E}}_e:=\overline{\tenb{U}}_e-\tena{1}_3=   \big(\tena \varphi,_j\cdot\, \tena d_i- \tena g_j\cdot \tena d_i^0\big)\,
  \tena d_i^0\dyad \tena g^j ,
\end{align*}
where $\tena{1}_3= \tena g_i\dyad \tena g^i= \tena d_i^0\dyad \tena d_i^0\,$ is the unit three-dimensional tensor.
As a strain measure for curvature (orientation change) one can employ the so-called \emph{wryness tensor} $\tenb \Gamma$ given by :
\begin{align}\label{e0,5}
    \tenb \Gamma:= \mathrm{axl}\big(\tenb Q_e^T\tenb Q_{e,i}\big)\dyad \tena g^i=  \tenb Q_0\big[  \mathrm{axl}\big(\overline{\tenb{R}}^T\overline{\tenb{R}}_{,i}\big) -  \mathrm{axl}\big(\tenb Q_0^T\tenb Q_{0,i}\big) \big] \dyad \tena g^i,
\end{align}
where $\mathrm{axl}\big(\tenb A\big)$ denotes the axial vector of any skew-symmetric tensor $\tenb A$.
For a detailed discussion on various strain measures of nonlinear micropolar continua we refer to the paper \cite{Pietraszkiewicz09}.

As an alternative to the wryness tensor $\,\tenb{\Gamma}\,$  one can make use of the $\,\mathrm{Curl}\,$ operator to define the so-called \emph{dislocation density tensor} $\,\overline{\tenb{D}}_e\,$ by \cite{Neff_curl08}
\begin{align}\label{e0,6}
    \overline{\tenb{D}}_e:=\tenb Q_e^T\,\mathrm{Curl}\,\tenb Q_e\,,
\end{align}
which is another curvature measure for micropolar continua. Note that the $\,\mathrm{Curl}\,$ operator has various definitions in the literature, but we will make its significance clear in the next Subsection 2.1, where we present the $\,\mathrm{Curl}\,$ operator in curvilinear coordinates.
The use of the dislocation density tensor $\,\overline{\tenb{D}}_e\,$ instead of the wryness tensor in conjuction with  micropolar and micromorphic media
has several advantages, as it was shown in \cite{GhibaNeffExistence,NeffGhibaMicroModel,MadeoNeffGhibaW}. The relationship between the wryness tensor $\,\tenb{\Gamma}\,$ and the dislocation density tensor $\,\overline{\tenb{D}}_e\,$ is discussed in the Subsection 2.2 in details.

Using the strain and curvature tensors $(\overline{\tenb{E}}_e\,,\,\overline{\tenb{D}}_e)$ the elastically stored energy density $W$ for the isotropic nonlinear Cosserat model can be expressed as \cite{Birsan-Neff-Ost_L56-2015,Lankeit-Neff-Ost-15}
\begin{alignat}{3}
      W(\overline{\tenb{E}}_e,\overline{\tenb{D}}_e) &= W_{\mathrm{mp}}(\overline{\tenb{E}}_e )+ W_{\mathrm{curv}}( \overline{\tenb{D}}_e),\qquad\mathrm{where} \vspace{4pt}\notag\\
      W_{\mathrm{mp}}(\overline{\tenb{E}}_e ) &= \mu\,\|\,\mathrm{dev_3\,sym}\, \overline{\tenb{E}}_e\,\|^2\, +  \,\mu_c \, \|\,\mathrm{skew}\, \overline{\tenb{E}}_e\,\|^2\,   + \, \dfrac{\kappa}{2}\,\big(\mathrm{tr}\,\overline{\tenb{E}}_e\big)^2\,,\vspace{4pt}\label{e0,7}\\
       W_{\mathrm{curv}}( \overline{\tenb{D}}_e) &= \mu\,L_c^p\,\Big(a_1\|\,\mathrm{dev_3\,sym}\, \overline{\tenb{D}}_e\|^2\, +  \,a_2 \|\,\mathrm{skew}\, \overline{\tenb{D}}_e\|^2\,   + \, a_3\big(\mathrm{tr}\,\overline{\tenb{D}}_e\big)^2\Big)^{{p}/{2}}\,,\notag
\end{alignat}
where $\mu$ is the shear modulus, $\kappa$ is the bulk modulus of classical isotropic elasticity, and  $\,\mu_c$ is called the \emph{Cosserat couple modulus}, which are assumed to satisfy
\begin{align*}
    \mu>0,\qquad \kappa>0,\qquad\text{and}\qquad \mu_c>0\,.
\end{align*}
The parameter $\,L_c$ introduces an internal length which is characteristic for the material, $a_i>0\,$ are dimensionless constitutive coefficients and $\,p\geq 2\,$ is a constant exponent.
Here, $\,\, \mathrm{dev}_3\,\tenb{X}:=\tenb{X}-\frac{1}{3}\,( \mathrm{tr}\,\tenb{X})\,\tenb{1}_3\,$ is the deviatoric part of any second order tensor $\tenb{X}$.

Under these assumptions on the constitutive coefficients the existence of minimizers to the corresponding minimization problem of the total energy functional has been shown, e.g., in \cite{Birsan-Neff-Ost_L56-2015,Lankeit-Neff-Ost-15}.

\subsection{The Curl operator}

For a vector field $\tena v$, the (coordinate-free) definition of the vector $\mathrm{curl}\, \tena v$ is
\begin{align}\label{e1}
  (\mathrm{curl}\, \tena v)\cdot \tena c= \mathrm{div}(\tena v\times \tena c)\qquad\text{for all constant vectors}\,\,\tena c,
\end{align}
where $\cdot$ denotes the scalar product and $\times$ the vector product. The $\mathrm{Curl}$ of a tensor field $\,\tenb T$ is the tensor field defined by
\begin{align}\label{e2}
  \big(\mathrm{Curl}\, \tenb T\big)^T \tena c= \mathrm{curl}\big(\tenb T^T \tena c\big)\qquad\text{for all constant vectors}\,\,\tena c.
\end{align}
\begin{remark}\label{rem1}
The operator $\mathrm{Curl}\,\tenb T$ given by \eqref{e2} coincides with the $\mathrm{Curl}$ operator defined in \cite{Svendsen02,Mielke06}. However, for other  authors the $\mathrm{Curl}$ of   $\,\tenb T$ is the transpose of $\mathrm{Curl}\,\tenb T$ defined by \eqref{e2}, see e.g.  \cite{Gurtin81,Erem-Leb-Alt13}.
\end{remark}

Then, from \eqref{e1} and \eqref{e2} we obtain the following formulas
\begin{align}\label{e3}
  \mathrm{curl}\, \tena v = -\tena v_{,i}\times \tena g^i,\qquad \mathrm{Curl}\, \tenb T= -\tenb T_{,i}\times \tena g^i.
\end{align}
Indeed, the definition \eqref{e1} yields
\begin{align*}
(\mathrm{curl}\, \tena v)\cdot \tena c  =\mathrm{div}(\tena v\times \tena c) = (\tena v\times \tena c)_{,i}\cdot \tena g^i=  (\tena v_{,i}\times \tena c)\cdot \tena g^i = ( \tena g^i \times \tena v_{,i})\cdot \tena c,
\end{align*}
and the equation \eqref{e3}$_1$ holds. Further, from \eqref{e2} we get
\begin{align*}
\big(\mathrm{Curl}\, \tenb T\big)^T \tena c= \mathrm{curl}\big(\tenb T^T \tena c\big) = \tena g^i \times\big(\tenb T^T \tena c\big)_{,i}=  \tena g^i \times\big(\tenb T_{,i}^T \tena c\big)=  \big(\tena g^i \times\tenb T_{,i}^T \big)\tena c,
\end{align*}
so it follows $\, \mathrm{Curl}\, \tenb T= \big(\tena g^i \times\tenb T_{,i}^T \big)^T= - \tenb T_{,i}\times \tena g^i$ and  the relations \eqref{e3} are proved.

In order to write the components of $\mathrm{curl}\, \tena v$ and $\mathrm{Curl}\,\tenb T$  in curvilinear coordinates, we introduce the following notations
\begin{align*}
g_{ij}=\tena g_i \dyad \tena g_j\,,\qquad g=\det\big( g_{ij}\big)_{3\times 3}\,>0.
\end{align*}
The alternating (Ricci) third-order tensor is
\begin{align*}
\tenb \epsilon=-\tenb 1_3\times \tenb 1_3= \epsilon_{ijk}\,\tena g^i \dyad \tena g^j\dyad \tena g^k= \epsilon^{ijk}\tena g_i \dyad \tena g_j\dyad \tena g_k\,,\qquad \mathrm{where}\\
\epsilon_{ijk}= \sqrt{g}\,e_{ijk}\,,\quad \epsilon^{ijk}= \dfrac{1}{\sqrt{g}}\,e_{ijk}\,,\quad e_{ijk}=\left\{ \begin{array}{rl}
                     1, & \,\,(i,j,k) \text{ is even permutation} \\
                 -1, & \,\,(i,j,k) \text{ is odd permutation} \\
              0, & \,\,(i,j,k) \text{ is no permutation}
         \end{array}
  \right. .
\end{align*}
The covariant, contravariant, and mixed components of any vector field $\tena v$ and any tensor field $\,\tenb T$ are introduced by
\begin{align*}
\tena v=v_k \,\tena g^k=v^k \tena g_k\,,\qquad \tenb T= T_{jk}\,\tena g^j\dyad \tena g^k= T^{jk}\tena g_j\dyad \tena g_k= T^{j}_{\cdot\, k} \,\tena g_j\dyad \tena g^k.
\end{align*}
For the partial derivatives with respect to $x_i$ we have the well-known expressions
\begin{align}\label{e4}
   \tena v_{,i} = v_{k|i} \,\tena g^k\,,\qquad \tenb T_{,i}= T_{jk|i}\,\tena g^j\dyad \tena g^k= T^{j}_{\cdot\, k|i} \,\tena g_j\dyad \tena g^k,
\end{align}
where a subscript bar preceding the index $i$ denotes covariant derivative w.r.t. $x_i$\,.

Using the relations \eqref{e4} in \eqref{e3} we can write the components of $\mathrm{curl}\, \tena v$ and $\mathrm{Curl}\,\tenb T$ as follows
\begin{align}\label{e5}
  \mathrm{curl}\, \tena v =\epsilon^{ijk}v_{j|i} \,\tena g_k\,,\qquad \mathrm{Curl}\, \tenb T = \epsilon^{ijk}T_{sj|i} \,\tena g^s\dyad \tena g_k = \epsilon^{ijk}T^{s}_{\cdot\,j|i} \,\tena g_s\dyad \tena g_k\,.
\end{align}
Indeed, from \eqref{e3}$_1$ and \eqref{e4}$_1$ we find
\begin{align*}
\mathrm{curl}\, \tena v = - \big( v_{k|i} \,\tena g^k \big)  \times \tena g^i = - v_{k|i} \big( \tena g^k  \times \tena g^i \big)  = - v_{k|i} \big( \epsilon^{kij}  \tena g_j \big)  = \epsilon^{ijk}v_{j|i}\,\tena g_k\,.
\end{align*}
Analogously, from \eqref{e3}$_2$ and \eqref{e4}$_2$ we get
\begin{align*}
\mathrm{Curl}\, \tenb T = -\big( T_{sk|i}\,\tena g^s\dyad \tena g^k  \big) \times \tena g^i = -T_{sk|i}\,\tena g^s\dyad \big( \tena g^k   \times \tena g^i\big) = \epsilon^{ijk}T_{sj|i} \,\tena g^s\dyad \tena g_k \,.
\end{align*}
Thus the equations \eqref{e5} are proved.
\begin{remark}\label{rem2}
In the special case of Cartesian coordinates, the relations \eqref{e3} and \eqref{e5} admit the simple form
\begin{align*}
\mathrm{curl}\, \tena v = -\tena v_{,i}\times \tena e_i= e_{ijk}v_{j,i} \,\tena e_k\,,\quad \mathrm{Curl}\, \tenb T= -\tenb T_{,i}\times \tena e_i= e_{ijk}\,T_{sj,i} \,\tena e_s\dyad \tena e_k \,,
\end{align*}
where $\tena v=v_i\tena e_i$ and   $\,\tenb T= T_{ij}\tena e_i\dyad \tena e_j$ are the corresponding coordinates. Moreover, in this case one can write
\begin{align}\label{e6}
   \mathrm{Curl}\,\tenb{T}  \,=\,  \tena{e}_i \dyad  \mathrm{curl}\big(\tenb{T}_i\big)\qquad\mathrm{for}\quad \tenb{T}= \tena{e}_i \dyad  \tenb{T}_i\,,
\end{align}
where $\, \tenb{T}_i=T_{ij}\,\tena{e}_j$ are the three  rows of the $3\times 3$ matrix $\,\big( T_{ij}\big)_{3\times 3}\,$.
The relation \eqref{e6} shows that $\,\mathrm{Curl}\,$ is defined row-wise  \cite{Neff_curl08}: the rows of the $3\times 3$ matrix $\,\mathrm{Curl}\,\tenb{T}\,$ are respectively the three vectors $\,\mathrm{curl}\big(\tenb{T}_i\big)$, \, $i=1,2,3$.
\end{remark}
\begin{remark}\label{rem3}
In order to write the corresponding formula in curvilinear coordinates which is analogous to \eqref{e6}, we introduce the vectors $\, \tenb{T}_i:=T_{ij}\,\tena{g}^j\,$ and $\, \tenb{T}^i:=T^{ij}\,\tena{g}_j\,=T^i_{\cdot\, j}\,\tena{g}^j \,$ such that it holds
\begin{align}\label{e7}
  \tenb{T}= \tena{g}^i \dyad  \tenb{T}_i\qquad\mathrm{and}\qquad  \tenb{T}= \tena{g}_i \dyad  \tenb{T}^i.
\end{align}
If we differentiate \eqref{e7}$_1$ with respect to $x_j$ we get
\begin{align*}
\tenb{T}_{,j} = \tena{g}^r_{\cdot\,,j} \dyad  \tenb{T}_r + \tena{g}^i \dyad  \tenb{T}_{i,j} = -\Gamma^r_{ji}\tena{g}^i  \dyad  \tenb{T}_r + \tena{g}^i \dyad  \tenb{T}_{i,j} = \tena{g}^i  \dyad \big( \tenb{T}_{i,j} -\Gamma^r_{ji}\,\tenb{T}_r  \big) ,
\end{align*}
where $\Gamma^r_{ij}$ are the Christoffel symbols of the second kind. Hence, it follows
\begin{align}\label{e8}
  \tenb{T}_{,j} =  \tena{g}^i \dyad  \tenb{T}_{i|j}\qquad\mathrm{with}\qquad  \tenb{T}_{i|j}:=  \tenb{T}_{i,j} -\Gamma^r_{ji}\,\tenb{T}_r = {T}_{ik|j} \,\tena{g}^k.
\end{align}
Taking the vector product of \eqref{e8}$_1$ with $\tena{g}^j$ we obtain
\begin{align*}
\mathrm{Curl}\, \tenb T = - \tenb T_{,j} \times \tena g^j =  -  \big( \tena{g}^i \dyad  \tenb{T}_{i|j}\big) \times \tena g^j,\quad\text{i.e.}
\end{align*}\vspace{-20pt}
\begin{align}\label{e9}
  \mathrm{Curl}\, \tenb T = \tena{g}^i \dyad  \mathrm{curl}_{\mathrm{cov}}\big(\tenb{T}_i\big)\qquad\mathrm{where}\quad \mathrm{curl}_{\mathrm{cov}}\big(\tenb{T}_i\big):= - \tenb{T}_{i|j}  \times \tena g^j.
\end{align}
The relation \eqref{e9} is the analogue of \eqref{e6} for curvilinear coordinates. Similarly, by differentiating \eqref{e7}$_2$ with respect to $x_j$ one can obtain the relation
\begin{align}\label{e10}
   \mathrm{Curl}\, \tenb T = \tena{g}_i \dyad  \mathrm{curl}_{\mathrm{cov}}\big(\tenb{T}^i\big)\qquad\text{where we denote}
\end{align}\vspace{-20pt}
\begin{align*}
   \mathrm{curl}_{\mathrm{cov}}\big(\tenb{T}^i\big):=  - \tenb{T}^{i}_{\cdot\,|j}  \times \tena g^j \qquad\mathrm{and}\quad  \tenb{T}^{i}_{\cdot\,|j}:= \tenb{T}^{i}_{\cdot\,,j} +  \Gamma^i_{rj}\,\tenb{T}^r = {T}^{i}_{\cdot\,k|j} \,\tena{g}^k.
\end{align*}
\end{remark}

\subsection{Relation between the   wryness tensor and the dislocation density tensor}

Let $\tenb A= A_{ij} \tena{g}^i \dyad\tena{g}^j$ be an arbitrary skew-symmetric tensor and $\mathrm{axl}(\tenb A)=a_k\tena{g}^k$ its axial vector. Then, the following relations hold
\begin{align}\label{e11}
\begin{array}{lcl}
  \tenb A & = & \mathrm{axl}(\tenb A)\times \tenb 1_3 =  \tenb 1_3  \times\mathrm{axl}(\tenb A), \vspace{4pt}\\
  \mathrm{axl}(\tenb A) & = & -\frac12 \,\tenb\epsilon : \tenb A = -\frac12 \, \epsilon^{ijk} A_{ij}\,\tena g_k\,,  \vspace{4pt}\\
  \tenb A & = &  - \tenb\epsilon\, \mathrm{axl}(\tenb A) =  -  \epsilon^{ijk} a_k\, \tena g_i\dyad \tena g_j\,,
\end{array}
\end{align}
where the double dot product `` : '' of two  tensors $\tenb{B}=B^{ijk}\,\tena{g}_i\dyad\tena{g}_j\dyad \tena{g}_k$ and $\tenb{T}=T_{ij}\,\tena{g}^i\dyad\tena{g}^j$
is defined as $\,\tenb{B}:\tenb{T}\,=\, B^{ijk}T_{jk}\,\tena{g}_i\,$.

Using these relations we can derive the close relationship between the   wryness tensor and the dislocation density tensor: it holds
\begin{align}\label{e12}
  \overline{\tenb{D}}_e = -\tenb \Gamma^T+(\mathrm{tr}\,\tenb \Gamma)\,\tenb 1_3\,\,,\quad \text{or equivalently},
\end{align}\vspace{-20pt}
\begin{align}\label{e13}
  \tenb \Gamma = - \overline{\tenb{D}}_e^T + \frac12 \, (\mathrm{tr}\,\overline{\tenb{D}}_e)\,\tenb 1_3\,\,.
\end{align}
Indeed, in view of the equation \eqref{e11}$_3$ and the definition \eqref{e0,5} we have
\begin{alignat*}{2}
\tenb Q_e^T\tenb Q_{e,k}\dyad \tena g^k & =  - \tenb\epsilon\, \mathrm{axl}(Q_e^T\tenb Q_{e,k})\dyad \tena g^k =  - \tenb\epsilon\,\tenb \Gamma \\
    & = - \big( \epsilon_{ijr}\,\tena g^i \dyad \tena g^j\dyad \tena g^r \big) \big( \Gamma^s_{\cdot\,k} \,\tena g_s \dyad \tena g^k\big) = -\epsilon_{ijs}\,\Gamma^s_{\cdot\,k} \,\tena g^i \dyad \tena g^j\dyad \tena g^k.
\end{alignat*}
Hence, we deduce
\begin{align}\label{e13,3}
\tenb Q_e^T\tenb Q_{e,k}\, =  -\epsilon_{ijs}\,\Gamma^s_{\cdot\,k} \,\tena g^i \dyad \tena g^j.
\end{align}
In view of \eqref{e3}$_2$\,, the definition \eqref{e0,6} can be written in the form
\begin{align}\label{e13,5}
   \overline{\tenb{D}}_e = \tenb Q_e^T\big(-\tenb Q_{e,k}\times \tena g^k\big) =  -\big(\tenb Q_e^T\tenb Q_{e,k}\big)\times \tena g^k .
\end{align}
Inserting \eqref{e13,3} in \eqref{e13,5}, we obtain
\begin{alignat*}{2}
\overline{\tenb{D}}_e & =  \epsilon_{ijs}\,\Gamma^s_{\cdot\,k} \,\big(\tena g^i \dyad \tena g^j\big)\times \tena g^k = \epsilon_{ijs}\,\Gamma^s_{\cdot\,k} \,\tena g^i \dyad \big( \epsilon^{jkr}\tena g_r \big) = \big(\epsilon_{jsi}\,\epsilon^{jkr} \big)\Gamma^s_{\cdot\,k} \,\tena g^i \dyad  \tena g_r \\
       & = \big( \delta^k_s \,\delta^r_i - \delta^r_s \,\delta^k_i   \big)\Gamma^s_{\cdot\,k} \,\tena g^i \dyad  \tena g_r = \Gamma^s_{\cdot\,s} \,\tena g^i \dyad  \tena g_i - \Gamma^s_{\cdot\,i} \,\tena g^i \dyad  \tena g_s = (\mathrm{tr}\,\tenb \Gamma)\,\tenb 1_3 -\tenb \Gamma^T.
\end{alignat*}
Thus, the relation \eqref{e12} is proved. If we apply the trace operator and the transpose in \eqref{e12} we obtain also the relation \eqref{e13}. For infinitesimal strains this formula is well-known under the name \emph{Nye's formula}, and $(\,-\,\tenb{\Gamma}\,)$ is also called \emph{Nye's curvature tensor} \cite{Nye53}. This relation has been first established in \cite{Neff_curl08}.

Let us find the components of the wryness tensor and the dislocation density tensor in curvilinear coordinates. To this aim, we write first the sqew-symmetric tensor
\begin{alignat}{2}\label{e14}
 \tenb Q_e^T\tenb Q_{e,i} & = \big(\tena d^0_j \dyad \tena d_j\big) \big(\tena d_{k,i} \dyad \tena d_k^0 + \tena d_{k} \dyad \tena d_{k,i}^0\big)= \big(  \tena d_j\cdot\tena d_{k,i} \big)\, \tena d^0_j \dyad  \tena d_k^0 + \tena d^0_j \dyad  \tena d_{j,i}^0 \notag\\
     & = \big(  \tena d_j\cdot\tena d_{k,i} - \tena d^0_j\cdot\tena d^0_{k,i} \big)\, \tena d^0_j \dyad  \tena d_k^0 \,.
\end{alignat}
Then, we obtain for the axial vector the equation
\begin{align}\label{e15}
  \mathrm{axl}\big(\tenb Q_e^T\tenb Q_{e,i}\big) = -\frac12 \, e_{jks}\, \big(  \tena d_j\cdot\tena d_{k,i} - \tena d^0_j\cdot\tena d^0_{k,i} \big)\, \tena d^0_s\,.
\end{align}
Indeed, according to \eqref{e11}$_2$ and \eqref{e14} we can write
\begin{alignat*}{3}
 \mathrm{axl}\big(\tenb Q_e^T\tenb Q_{e,i}\big) & =  -\frac12 \,\tenb\epsilon :  \big(\tenb Q_e^T\tenb Q_{e,i}\big) \\
      & =   -\frac12 \,\big(e_{sjk}\,\tena d_s^0 \dyad \tena d_j^0 \dyad \tena d_k^0 \big) :  \big[ \big(  \tena d_l\cdot\tena d_{r,i} - \tena d^0_l\cdot\tena d^0_{r,i} \big)\, \tena d^0_l \dyad  \tena d_r^0 \,\big] \\
       & =   -\frac12 \, e_{jks}\, \big(  \tena d_j\cdot\tena d_{k,i} - \tena d^0_j\cdot\tena d^0_{k,i} \big)\, \tena d^0_s
\end{alignat*}
and the relation \eqref{e15} is proved. Using \eqref{e15}  in the definition \eqref{e0,5}  we find the following formula for the wryness tensor
\begin{align}\label{e16}
  \tenb\Gamma\,=\, \frac12 \, e_{jks}\, \big(  \tena d_{j,i}\cdot\tena d_{k} - \tena d^0_{j,i}\cdot\tena d^0_{k} \big)\, \tena d^0_s\dyad \tena g^i.
\end{align}
To obtain an expression for the components of $\overline{\tenb{D}}_e$  we insert \eqref{e14} in  \eqref{e13,5}  and we get
\begin{alignat}{2}\label{e18}
 \overline{\tenb{D}}_e  & = - \big(  \tena d_j\cdot\tena d_{k,i} - \tena d^0_j\cdot\tena d^0_{k,i} \big)\, \big( \tena d^0_j \dyad  \tena d_k^0\big)\times \tena g^i \notag \\
        & = \big(  \tena d_{j,i}\cdot\tena d_{k} - \tena d^0_{j,i}\cdot\tena d^0_{k} \big)\, \tena d^0_j \dyad  \big(  \tena d_k^0\times \tena g^i \big).
\end{alignat}
We rewrite the last vector product as
\begin{align*}
\tena d_k^0\times \tena g^i = \tena d_k^0\times \big[ \big( \tena g^i \cdot \tena d_r^0 \big)\,\tena d_r^0  \big] = \big( \tena g^i \cdot \tena d_r^0 \big)\,\tena d_k^0\times   \tena d_r^0  =  e_{krs} \big( \tena g^i \cdot \tena d_r^0 \big)\,\tena d_s^0
\end{align*}
and we insert it in \eqref{e18} to find the following expression for the dislocation density tensor
\begin{align}\label{e19}
  \overline{\tenb{D}}_e =  e_{krs} \big(  \tena d_{j,i}\cdot\tena d_{k} - \tena d^0_{j,i}\cdot\tena d^0_{k} \big)\,   \big( \tena g^i \cdot \tena d_r^0 \big)\, \tena d^0_j \dyad\tena d_s^0\,.
\end{align}
\begin{remark}\label{rem4}
In the special case of Cartesian coordinates one can identify $\tena d_i^0=\tena e_i\,,\,\tena g^i=\tena g_i=\tena e_i\,$, and the relations \eqref{e16} and \eqref{e18} simplify to the forms
\begin{alignat*}{2}
\tenb\Gamma & =\, \frac12 \, e_{iks}\, \big(  \tena d_{k,j}\cdot\tena d_{s} \big)\, \tena e_i\dyad \tena e_j\,\,, \\
 \overline{\tenb{D}}_e & =  e_{ijk} \, \big(  \tena d_{j,i}\cdot\tena d_{s} \big)\, \tena e_s\dyad \tena e_k\,\,.
\end{alignat*}
\end{remark}
\begin{remark}\label{rem5}
One can find various definitions of the wryness tensor in the literature, see e.g. \cite{Tambaca10}, where $\tenb \Gamma$ is called the \emph{curvature strain tensor}. Thus, one can alternatively define the wryness tensor by
\begin{align}\label{e20}
 \tenb\Gamma = \tenb Q_e^T\,\tenb\omega\,,
\end{align}
where $\tenb \omega$ is the second order tensor given by
\begin{align}\label{e21}
  \tenb\omega = \tena \omega_i\dyad \tena g^i \qquad\text{with}\qquad \tenb Q_{e,i} = \tena\omega_i\times \tenb Q_e\,\,.
\end{align}
If we compare the definition \eqref{e0,5} with \eqref{e20}, \eqref{e21}, we see that indeed $\tenb Q_e^T\,\tena\omega_i= \mathrm{axl}\big(\tenb Q_e^T\tenb Q_{e,i}\big)$ , i.e.
\begin{align}\label{e22}
  \tena\omega_i\, = \,\tenb Q_e\,\mathrm{axl}\big(\tenb Q_e^T\tenb Q_{e,i}\big)   =    \mathrm{axl}\big(\tenb Q_{e,i}\,\tenb Q_e^T\big).
\end{align}
By a straightforward but lengthy calculation one can prove that the vectors $\tena \omega_i$ are expressed in terms of the directors by
\begin{align}\label{e23}
   \tena\omega_i\, = \,\frac12 \big[   \tena d_{j}\times  \tena d_{j,i} -  \tenb Q_e \big(\tena d^0_{j}\times  \tena d^0_{j,i} \big)  \big].
\end{align}
Inserting \eqref{e23} in \eqref{e21}$_1$ and \eqref{e20} we obtain the expression of the wryness tensor written with the help of the directors $\tena d_i$
\begin{align}\label{e24}
  \tenb\Gamma \, = \, \,\frac12 \big[  \tenb Q_e^T \big(\tena d_{j}\times  \tena d_{j,i}\big) -    \tena d^0_{j}\times  \tena d^0_{j,i}   \big]\dyad \tena g^i .
\end{align}
\end{remark}

\section{The Curl operator on surfaces}\label{sec:3}

Let $\mathcal{S}$ be a smooth surface embedded in the Euclidean space $\mathbb{R}^3$ and let $\tena y_0(x_1,x_2)$, $\tena y_0:\omega\rightarrow\mathbb{R}^3$, be a parametrization of this surface. We denote the covariant base vectors in the tangent plane by $\tena a_1\,,\tena a_2$ and the contravariant base vectors   by $\tena a^1\,,\tena a^2\,$:
\begin{align*}
 \tena a_\alpha= \dfrac{\partial \tena y_0}{\partial x_\alpha }\,=\tena y_{0,\alpha}\,,\qquad \tena a_\alpha\cdot \tena a^\beta=\delta_\alpha^\beta
\end{align*}
and let
\begin{align*}
 \tena a_3= \tena a^3= \tena n_0= \dfrac{\tena a_1\times\tena a_2}{|\tena a_1\times\tena a_2| }\,,
\end{align*}
where $\tena n_0$ is the unit normal to the surface. Further, we designate by
\begin{align*}
 a_{\alpha\beta}=\tena{a}_\alpha\cdot \tena{a}_\beta\,,\qquad a^{\alpha\beta}=\tena{a}^\alpha\cdot \tena{a}^\beta,\qquad   a =\sqrt{\mathrm{det}\, \big( a_{\alpha\beta} \big){}_{2\times 2}}\,= |\tena a_1\times\tena a_2|\,>0
\end{align*}
and we have
\begin{align}\label{e26}
   \tena{a}^\alpha\times \tena{a}^\beta= \epsilon^{\alpha\beta} \tena a_3\,,\,\,\, \tena{a}^3\times \tena{a}^\alpha= \epsilon^{\alpha\beta} \tena a_\beta\,,\,\,\, \tena{a}_\alpha\times \tena{a}_\beta= \epsilon_{\alpha\beta} \tena a^3,\quad \tena{a}_3\times \tena{a}_\alpha= \epsilon_{\alpha\beta} \tena a^\beta,
\end{align}
where $\epsilon^{\alpha\beta} =\,\dfrac{1}{a}\,e_{\alpha\beta} \,,\,\, \epsilon_{\alpha\beta} ={a}\,e_{\alpha\beta} $ and $\,e_{\alpha\beta} $ is the two-dimensional alternator given by $e_{12}=-e_{21}=1, \,e_{11}=e_{22}=0$.

Then, $\tena{a}=    a_{\alpha\beta}\tena{a}^\alpha\dyad \tena{a}^\beta= a^{\alpha\beta}\tena{a}_\alpha\dyad \tena{a}_\beta=\tena{a}_\alpha\dyad \tena{a}^\alpha$ represents the first fundamental tensor of the surface $\mathcal{S}$, while the second fundamental tensor $\tena b$ is defined by
\begin{alignat*}{2}
& \tena{b}= -\mathrm{Grad}_s\,\tena{n}_0=-  \tena{n}_{0,\alpha}\dyad\tena{a}^\alpha= b_{\alpha\beta}\,\tena{a}^\alpha\dyad \tena{a}^\beta=b^\alpha_\beta\,\tena{a}_\alpha\dyad \tena{a}^\beta, \quad\text{with}\vspace{4pt}\\
&    b_{\alpha\beta}=-  \tena{n}_{0,\beta}\cdot\tena{a}_\alpha= b_{\beta\alpha}\,,\qquad b^{\alpha}_{\beta}=-  \tena{n}_{0,\beta}\cdot\tena{a}^\alpha\,.
\end{alignat*}
The surface gradient $\mathrm{Grad}_s$ and surface divergence $\mathrm{Div}_s$ operators are defined for a vector field $\tena v$ by
\begin{align}\label{e27}
  \mathrm{Grad}_s\,\tena v\, =\, \dfrac{\partial\tena v}{\partial x_\alpha}\,\dyad \tena a^\alpha= \tena v_{,\alpha}\dyad \tena a^\alpha,\qquad \mathrm{Div}_s\,\tena v\, =\,\mathrm{tr}\big[ \mathrm{Grad}_s\,\tena v\big]= \tena v_{,\alpha}\cdot \tena a^\alpha.
\end{align}
We also introduce the so-called \emph{alternator tensor} $\tena c$ of the surface \cite{Zhilin06}
\begin{align}\label{e28}
   \tena c=-\tena n_0\times \tena a=-\tena a\times\tena n_0=   \epsilon^{\alpha\beta}\,\tena{a}_\alpha\dyad \tena{a}_\beta =  \epsilon_{\alpha\beta}\,\tena{a}^\alpha\dyad \tena{a}^\beta.
\end{align}
The tensors $\tena a$ and $\tena b$ are symmetric, while $\tena c$ is skew-symmetric and satisfies  $\, \tena c\tena c=-\tena a$. Note that the tensors $\tena a\,$, $\tena b\,$, and $\tena c$ defined above are \emph{planar}, i.e. they are tensors in the tangent plane of the surface. Moreover, $\tena a$ is the identity tensor in the tangent plane.\medskip

We define \textbf{the surface Curl operator} $\mathrm{curl}_s$ for vector fields $\tena v$ and, respectively, $\mathrm{Curl}_s$ for tensor fields $\tenb T$ by
\begin{alignat}{2}
& \big(\mathrm{curl}_s\,\tena v\big)\cdot \tena k \, = \,  \mathrm{Div}_s\big(\tena v\times \tena k\big)\qquad\text{for all constant vectors }\tena k, \label{e29}
\vspace{4pt}\\
&  \big(\mathrm{Curl}_s\,\tenb T\big)^T \tena k \, = \,  \mathrm{curl}_s\big(\tenb T^T \tena k\big)\qquad\text{for all constant vectors }\tena k. \label{e30}
\end{alignat}
Thus, $\mathrm{curl}_s\,\tena v$ is a vector field, while $\mathrm{Curl}_s\,\tenb T$ is a tensor field.
\begin{remark}\label{rem6}
These definitions are analogous to the corresponding definitions \eqref{e1}, \eqref{e2} in the three-dimensional case. Notice that the curl operator on surfaces has a different significance for other authors, see e.g. \cite{Book-Geomagn}.
\end{remark}
From the definitions \eqref{e29} and  \eqref{e30} it follows
\begin{align}\label{e31}
    \mathrm{curl}_s\,\tena v = - \tena v_{,\alpha}\times \tena a^\alpha,\qquad \mathrm{Curl}_s\,\tenb T = - \tenb T_{,\alpha}\times \tena a^\alpha.
\end{align}
Indeed, in view of \eqref{e27} and  \eqref{e29} we have
\begin{alignat*}{2}
\big(\mathrm{curl}_s\,\tena v\big)\cdot \tena k & =   \mathrm{Div}_s\big(\tena v\times \tena k\big) = \big(\tena v\times \tena k\big)_{,\alpha}\cdot \tena a^\alpha = \big(\tena v_{,\alpha}\times \tena k\big)\cdot \tena a^\alpha \\
& = \big(\tena a^\alpha\times \tena v_{,\alpha}\big)\cdot \tena k = \big(- \tena v_{,\alpha}\times \tena a^\alpha\big)\cdot \tena k\quad \text{for all constant vectors }\tena k
\end{alignat*}
and also
\begin{align*}
\big(\mathrm{Curl}_s\,\tenb T\big)^T \tena k  =   \mathrm{curl}_s\big(\tenb T^T \tena k\big) = \tena a^\alpha\times \big(\tenb T^T \tena k\big)_{,\alpha} = \tena a^\alpha\times \big(\tenb T^T_{,\alpha} \tena k\big) = \big(\tena a^\alpha\times \tenb T^T_{,\alpha}\big) \tena k,
\end{align*}
which implies $ \mathrm{Curl}_s\,\tenb T = \big(\tena a^\alpha\times \tenb T^T_{,\alpha}\big)^T = - \tenb T_{,\alpha}\times \tena a^\alpha $, so the relations \eqref{e31} hold true.

To write the components of $\mathrm{curl}_s\,\tena v$ and $\mathrm{Curl}_s\,\tenb T$ we employ the covariant derivatives on the surface. Let $\tena v= v_i\,\tena a^i$ be a vector field on $\mathcal{S}$. Then, we have
\begin{alignat}{2}
   & \tena a^\alpha_{\cdot ,\beta} = - \Gamma^\alpha_{\beta\gamma}\, \tena a^\gamma + b^\alpha_\beta\,\tena a^3,\qquad \tena a_{3,\beta} = - b^\alpha_\beta\,\tena a_\alpha= -b_{\alpha\beta}\,\tena a^\alpha,\notag\\
   &  \tena v_{,\alpha} = ( v_{\beta|\alpha} - b_{\alpha\beta}\,v_3) \tena a^\beta + (v_{3,\alpha}+ b_\alpha^\beta\,v_\beta)\tena a^3, \label{e31,5}
\end{alignat}
where $ v_{\beta|\alpha} = v_{\beta,\alpha} - \Gamma^\gamma_{\alpha\beta}\,v_\gamma\,$ is the covariant derivative with respect to $x_\alpha\,$. Inserting this relation in \eqref{e31}$_1$ and using \eqref{e26}$_{1,2}$ we obtain
\begin{align}\label{e32}
    \mathrm{curl}_s\,\tena v\,=\, \epsilon^{\alpha\beta}\big[ (v_{3,\beta}+ b^\gamma_\beta\,v_\gamma)\tena a_\alpha +  v_{\beta|\alpha}\, \tena a_3 \big].
\end{align}
For a tensor field $\tenb T=T_{ij}\,\tena a^i\dyad\tena a^j = T^{ij}\,\tena a_i\dyad\tena a_j = T^{i}_{\cdot\,j}\,\tena a_i\dyad\tena a^j$ on the surface, the derivative $\tenb T_{,\gamma}$ can be expressed as
\begin{align}\label{e33}
   \tenb T_{,\gamma} = \big( T_{\alpha\beta|\gamma} - b_{\alpha\gamma}\,T_{3\beta} - b_{\beta\gamma}\,T_{\alpha 3} \big) \tena a^\alpha\dyad\tena a^\beta + \big( T_{\alpha 3|\gamma} + b^\beta_{\gamma}\,T_{\alpha\beta} - b_{\alpha\gamma}\,T_{3 3} \big) \tena a^\alpha\dyad\tena a^3\notag\\
   + \big( T_{3\alpha|\gamma}+ b^\beta_{\gamma}\,T_{\beta\alpha} - b_{\alpha\gamma}\,T_{3 3} \big) \tena a^3\dyad\tena a^\alpha + \big( T_{3 3,\gamma} + b^\alpha_{\gamma}\,T_{\alpha 3} +  b^\alpha_{\gamma}\,T_{3\alpha } \big) \tena a^3\dyad\tena a^3,
\end{align}
where the covariant derivatives are
\begin{alignat*}{2}
& T_{\alpha\beta|\gamma} =  T_{\alpha\beta,\gamma} - \Gamma^\delta_{\beta\gamma}\, T_{\alpha\delta} - \Gamma^\delta_{\alpha\gamma}\, T_{\delta\beta}\,\,,\\
  & T_{\alpha 3|\gamma} =  T_{\alpha 3,\gamma} - \Gamma^\beta_{\alpha\gamma}\, T_{\beta 3}\,\,,\qquad T_{3\alpha |\gamma} =  T_{3\alpha ,\gamma} - \Gamma^\beta_{\alpha\gamma}\, T_{3\beta }\,\,.
\end{alignat*}
Using \eqref{e33}  in \eqref{e31}$_2$ we obtain with the help of \eqref{e26}$_{1,2}$
\begin{align}\label{e34}
    \mathrm{Curl}_s\,\tenb T = \epsilon^{\beta\gamma}\big(T_{\alpha 3|\gamma} \!+\! b^\sigma_{\gamma}\,T_{\alpha\sigma} \!-\! b_{\alpha\gamma}\,T_{3 3} \big) \tena a^\alpha\!\dyad\!\tena a_\beta  + \epsilon^{\gamma\beta}\big(T_{\alpha\beta|\gamma}\! - \!b_{\alpha\gamma}\,T_{3\beta} \big) \tena a^\alpha\!\dyad\!\tena a_3 \notag\\
    +  \epsilon^{\beta\gamma}\big(T_{3 3,\gamma} + b^\alpha_{\gamma}\,T_{\alpha 3} +  b^\alpha_{\gamma}\,T_{3\alpha } \big) \tena a^3\!\dyad\!\tena a_\beta + \epsilon^{\gamma\beta}\big(T_{3\beta|\gamma}\! + \!b^\alpha_{\gamma}\,T_{\alpha\beta} \big) \tena a^3\!\dyad\tena a_3\,.
\end{align}
Alternatively, one can use the mixed components $T^{i}_{\cdot\,j}\,$ and write $\mathrm{Curl}_s\,\tenb T$ in the tensor basis $\{ \,\tena a_i\dyad\tena a_j  \}$
\begin{align}\label{e35}
    \mathrm{Curl}_s\,\tenb T = \epsilon^{\beta\gamma}\big(T^{\alpha}_{\,\cdot\, 3|\gamma} \!+\! b^\sigma_{\gamma}\,T^{\alpha}_{\,\cdot\, \sigma} \!-\! b^\alpha_{\gamma}\,T^3_{\,\cdot\, 3} \big) \tena a_\alpha\!\dyad\!\tena a_\beta  + \epsilon^{\gamma\beta}\big(T^{\alpha}_{\,\cdot\, \beta|\gamma}\! - \!b^\alpha_{\gamma}\,T^3_{\,\cdot\,\beta} \big) \tena a_\alpha\!\dyad\!\tena a_3 \notag\\
    +  \epsilon^{\beta\gamma}\big(T^3_{\,\cdot\, 3,\gamma} + b_{\alpha\gamma}\,T^{\alpha}_{\,\cdot\,  3} +  b^\alpha_{\gamma}\,T^3_{\,\cdot\,\alpha } \big) \tena a_3\!\dyad\!\tena a_\beta + \epsilon^{\gamma\beta}\big(T^3_{\,\cdot\,\beta|\gamma}\! + \!b_{\alpha\gamma}\,T^{\alpha}_{\,\cdot\, \beta} \big) \tena a_3\!\dyad\tena a_3\,.
\end{align}
where
\begin{alignat*}{2}
& T^{\alpha}_{\,\cdot\, \beta|\gamma} =  T^{\alpha}_{\,\cdot\, \beta,\gamma} + \Gamma^\alpha_{ \gamma\sigma}\, T^{\sigma}_{\,\cdot\, \beta} - \Gamma^\sigma_{\beta\gamma}\, T^{\alpha}_{\,\cdot\, \sigma}\,\,,\\
  & T^{\alpha}_{\,\cdot\,  3|\gamma} =  T^{\alpha}_{\,\cdot\,  3,\gamma} + \Gamma^\alpha_{\gamma\sigma}\, T^{\sigma}_{\,\cdot\,  3}\,\,,\qquad T^{3}_{\,\cdot\, \beta |\gamma} =  T^{3}_{\,\cdot\, \beta ,\gamma} - \Gamma^\sigma_{\beta\gamma}\, T^{3}_{\,\cdot\, \sigma }\,\,.
\end{alignat*}
\begin{remark}\label{rem7}
In order to obtain a formula analogous to \eqref{e6} and \eqref{e9}, \eqref{e10} for $ \mathrm{Curl}_s\,$ we write $\tenb T$ in the form
\begin{align*}
\tenb T =   \tena a^i\dyad   \tenb T_i =  \tena a_i\dyad   \tenb T^i     \qquad \text{with}  \qquad                \tenb T_i = T_{ij}\, \tena a^j ,\quad  \tenb T^i = T^{i}_{\,\cdot j}\,\tena a^j.
\end{align*}
By differentiating the first equation with respect to $x_\gamma$ we get
\begin{alignat*}{2}
\tenb T_{,\gamma} & =   \tena a^i_{\,\,,  \gamma}\!\dyad \!  \tenb T_i +  \tena a^i\!\dyad\!   \tenb T_{i,  \gamma} = \big(\! -\!\Gamma^\alpha_{\beta\gamma}\,\tena a^\beta +  b^\alpha_{\gamma}\,\tena a^3\big) \!\dyad \!  \tenb T_\alpha -b_{\alpha\gamma}\,\tena a^\alpha \!\dyad\!   \tenb T_3 + \tena a^i\!\dyad \!  \tenb T_{i,  \gamma} \\
& = \tena a^\alpha\dyad \big(\tenb T_{\alpha,  \gamma} -\Gamma^\beta_{\alpha\gamma}\,\tenb T_\beta -  b_{\alpha\gamma}\,\tenb T_3\big) + \tena a^3\dyad \big(\tenb T_{3,  \gamma} + b^\alpha_{\gamma}\,\tenb T_\alpha\big).
\end{alignat*}
Taking the vector product with $\tena a^\gamma$ and using \eqref{e31}$_2$ we find
\begin{align}\label{e36}
    \mathrm{Curl}_s\,\tenb T = -  \big[  \tena a^\alpha\dyad \big(\tenb T_{\alpha|  \gamma} -  b_{\alpha\gamma}\,\tenb T_3\big) + \tena a^3\dyad \big(\tenb T_{3,  \gamma} + b^\alpha_{\gamma}\,\tenb T_\alpha\big)  \big]\times\tena a^\gamma,
\end{align}
with $\,\tenb T_{\alpha|  \gamma}\,:=\, \tenb T_{\alpha,  \gamma} -\Gamma^\beta_{\alpha\gamma}\,\tenb T_\beta\,$. Similarly, we obtain
\begin{align}\label{e37}
    \mathrm{Curl}_s\,\tenb T = -  \big[  \tena a_\alpha\dyad \big(\tenb T^{\alpha}_{\,\cdot\, |  \gamma} -  b^\alpha_{\gamma}\,\tenb T^3\big) + \tena a_3\dyad \big(\tenb T^{3}_{\,\cdot\, ,  \gamma} + b_{\alpha\gamma}\,\tenb T^\alpha\big)  \big]\times\tena a^\gamma,
\end{align}
with $\,\tenb T^{\alpha}_{\,\cdot\, |  \gamma}\,:=\, \tenb T^{\alpha}_{\,\cdot\, ,  \gamma} + \Gamma^\alpha_{\beta\gamma}\,\tenb T^\beta\,$. The equations \eqref{e36} and \eqref{e37} are the counterpart of the relations \eqref{e9} and, respectively, \eqref{e10} in the three-dimensional theory.
\end{remark}

\section{The shell dislocation density tensor}\label{sec:4}

Let us present first the kinematics of Cosserat-type shells, which coincides with the kinematics of the 6-parameter shell model, see \cite{Pietraszkiewicz-book04,Eremeyev06,Birsan-Neff-L54-2014}.

We consider a deformable surface $\omega_\xi\subset\mathbb{R}^3$ which is identified with the midsurface of the shell in its reference configuration and denote with $(\xi_1,\xi_2,\xi_3)$ a generic point of the surface. Each material point is assumed to have 6 degrees of freedom (3 for translations and 3 for rotations). Thus, the deformation of the Cosserat-type shell is determined by a vectorial map $\tena m_\xi$ and the microrotation tensor $\tenb{R}_\xi\,$
\begin{align*}
\tena m_\xi:\omega_\xi \rightarrow\omega_c\,, \qquad \tenb{R}_\xi:\omega_\xi \rightarrow \mathrm{SO}(3)\, ,
\end{align*}
where $\omega_c$ denotes the deformed (current) configuration of the surface. We consider a parametric representation $\tena y_0$ of the reference configuration $\omega_\xi$
\begin{align*}
\tena y_0 :\omega \rightarrow\omega_\xi\,, \qquad \tena y_0(x_1,x_2)=(\xi_1,\xi_2,\xi_3),
\end{align*}
where $\omega\subset\mathbb{R}^2$ is the bounded variation domain (with Lipschitz boundary) of the parameters $(x_1,x_2)$. Using the same notations as in Section \ref{sec:3}, we introduce the base vectors $\tena a_i\,,\,\tena a^j$ and the fundamental tensors $\tena a\,,\,\tena b$ for the reference surface $\omega_\xi\,$.

The deformation function $\tena m$ is then defined by the composition
\begin{align*}
\tena m= \tena m_\xi\circ\tena y_0 :\omega \rightarrow\omega_c\,,\qquad \tena m(x_1,x_2): = \tena m_\xi\big(\tena y_0(x_1,x_2)\big).
\end{align*}
According to \eqref{e27}, the surface gradient of the deformation has the expression
\begin{align}\label{e38}
      \mathrm{Grad}_s\,\tena m\, =\,   \tena m_{,\alpha}\dyad \tena a^\alpha.
\end{align}
As in the three-dimensional case (see Section \ref{sec:2}) we define the \emph{elastic microrotation} $\tenb Q_e$ by the composition
\begin{align*}
\tenb Q_e= \tenb R_\xi\circ\tena y_0 :\omega \rightarrow \mathrm{SO}(3 ),\qquad \tenb Q_e (  x_1,x_2):= \tenb R_\xi\big(\tena y_0( x_1,x_2)\big),
\end{align*}
the \emph{total microrotation} $\overline{\tenb R}$ by
\begin{align*}
\overline{\tenb R} :\omega \rightarrow \mathrm{SO}(3 ),\qquad \overline{\tenb R} (x_1,x_2)= \tenb Q_e (  x_1,x_2)\,\tenb Q_0 (  x_1,x_2),
\end{align*}
where $\tenb Q_0  :\omega \rightarrow \mathrm{SO}(3 )$  is the \emph{initial microrotation}, which describes the orientation of points in the reference configuration.

To characterize  the orientation and rotation of points in Cosserat-type shells one  employs (as in the three-dimensional case) a triad of orthonormal directors attached to each point. We denote by $\tena d_i^0(x_1,x_2) $ the directors in the reference configuration $\omega_\xi$  and by $\tena d_i(x_1,x_2) $ the directors in the deformed configuration $\omega_c\,$($i=1,2,3$). The domain $\omega$ is refered to an orthogonal Cartesian frame $Ox_1x_2x_3$ such that $\omega\subset Ox_1x_2\,$ and let $\tena e_i\,$ be the unit vectors along the coordinate axes $Ox_i\,$. Then, the microrotation tensors can be expressed as follows
\begin{align}\label{e39}
    \tenb Q_e=\tena d_i\dyad\tena d_i^0\,,\qquad \overline{\tenb R} = \tenb Q_e \,\tenb Q_0 = \tena d_i\dyad\tena e_i\,,\qquad \tenb Q_0 = \tena d_i^0\dyad\tena e_i\,.
\end{align}
\begin{remark}\label{rem8}
The initial directors $\tena d_i^0$ are usually chosen such that
\begin{align}\label{e39,5}
     \tena d_3^0 = \tena n_0  \,,\qquad \tena d_\alpha^0\cdot \tena n_0 = 0\,,
\end{align}
i.e. $\tena d_3^0$ is orthogonal to $\omega_\xi$ and $\tena d_\alpha^0$ belong to the tangent plane. This assumption is not necessary in general, but it will be adopted here since it simplifies many of the subsequent expressions. In the deformed configuration, the director $\tena d_3\,$ is no longer orthogonal to the surface $\omega_c\,$ (the Kirchhof-Love condition is not imposed). One convenient  choice of the initial microrotation tensor $\tenb Q_0 = \tena d_i\dyad\tena e_i\,$ such that the conditions \eqref{e39,5} be satisfied is
$\tenb Q_0 = \mathrm{polar}\big(\tena a_i\dyad\tena e_i\big)$, as it was shown in Remark 10 of \cite{Birsan-Neff-MMS-2014}.
\end{remark}

Let us present next the shell strain and curvature measures. In the 6-parameter shell theory the \emph{elastic shell strain tensor} $\tenb E_e$ is defined by \cite{Pietraszkiewicz-book04,Eremeyev06}
\begin{align}\label{e40}
      \tenb E_e\,=\, \tenb Q_e^T\, \mathrm{Grad}_s\tena m - \tena a\,.
\end{align}
To write the components of $\tenb E_e$ we insert \eqref{e38} and \eqref{e39}$_1$ into \eqref{e40}
\begin{align*}
\tenb E_e\,=\, \big( \tena d_i^0\dyad\tena d_i  \big)\big( \tena m_{,\alpha}\dyad \tena a^\alpha  \big)  - \tena a_{\alpha}\dyad \tena a^\alpha    =  \big( \tena m_{,\alpha}\cdot\tena d_i- \tena a_\alpha\cdot\tena d_i^0   \big)\, \tena d_i^0\dyad \tena a^\alpha.
\end{align*}
As a measure of orientation (curvature) change, the \emph{elastic shell bending-curvature tensor} $\tenb K_e$ is defined by \cite{Pietraszkiewicz-book04,Eremeyev06,Birsan-Neff-L54-2014}
\begin{align}\label{e41}
    \tenb K_e =  \mathrm{axl}\big( \tenb Q_e^T\tenb Q_{e,\alpha}  \big) \dyad \tena a^\alpha = \tenb Q_0 \big[ \mathrm{axl}\big( \overline{\tenb R}^T\overline{\tenb R}_{,\alpha}  \big) - \mathrm{axl}\big( \tenb Q_0^T\tenb Q_{0,\alpha}  \big)   \big].
\end{align}
We remark the analogy to the definition \eqref{e0,5} of the wryness tensor $\tenb \Gamma$ in the three-dimensional theory. Following the analogy to \eqref{e0,6}, we employ next the surface curl operator $\mathrm{Curl}_s$  defined in Section \ref{sec:3} to introduce the new \emph{shell dislocation density tensor} $\tenb D_e$ by
\begin{align}\label{e42}
      \tenb D_e = \tenb Q_e^T\,\mathrm{Curl}_s\,\tenb Q_e\,\,.
\end{align}
In view of relation \eqref{e31}$_2\,$, we can write this definition in the form
\begin{align}\label{e43}
      \tenb D_e = \tenb Q_e^T \big(-\tenb Q_{e,\alpha}\times \tena a^\alpha \big)=  -\big(\tenb Q_e^T\tenb Q_{e,\alpha}\big)\times \tena a^\alpha .
\end{align}
The tensor $\tenb D_e$ given by \eqref{e42} represents an alternative strain measure for orientation (curvature) change in Cosserat-type shells.

In what follows, we want to establish the relationship between the shell bending-curvature tensor $\tenb K_e$ and the shell dislocation density tensor $\tenb D_e\,$. We observe that this relationship is analogous to the corresponding relations \eqref{e14}, \eqref{e15} in the three-dimensional theory. More precisely, in the shell theory it holds
\begin{align}\label{e44}
      \tenb D_e = - \tenb K_e^T+ \big(\mathrm{tr}\,\tenb K_e\big)\tena 1_3\qquad\!\text{or equivalently,}\!\qquad \tenb K_e = - \tenb D_e^T+ \frac12\big(\mathrm{tr}\,\tenb D_e\big)\tena 1_3.
\end{align}
To prove \eqref{e44}, we designate the components of the shell bending-curvature tensor  by $\tenb K_e = K_{i\alpha}\,\tena d_i^0\dyad\tena a^\alpha $  and use \eqref{e13}$_3$ to write
\begin{alignat*}{2}
\big(\tenb Q_e^T\tenb Q_{e,\alpha}\big)&\dyad \tena a^\alpha  = - \tenb \epsilon\,  \mathrm{axl}\big( \tenb Q_e^T\tenb Q_{e,\alpha}  \big) \dyad \tena a^\alpha = - \tenb \epsilon\, \tenb K_e \\
 & = - \big( e_{ijk}\, \tena d_i^0 \dyad \tena d_j^0 \dyad \tena d_k^0    \big) \big(  K_{s\alpha}\,\tena d_s^0\dyad\tena a^\alpha  \big) =  -   e_{ijs}\,     K_{s\alpha}\,\tena d_i^0 \dyad \tena d_j^0 \dyad\tena a^\alpha,
\end{alignat*}
which implies
\begin{align*}
\tenb Q_e^T\tenb Q_{e,\alpha}\,=\,  -   e_{ijs}\,     K_{s\alpha}\,\tena d_i^0 \dyad \tena d_j^0
\end{align*}
We substitute the last relation into \eqref{e43} and derive
\begin{alignat*}{6}
\tenb D_e & =  \big( e_{ijs}\, K_{s\alpha}\,\tena d_i^0 \dyad \tena d_j^0 \big)\times \tena a^\alpha = \big( e_{ijs}\, K_{s\alpha}\,\tena d_i^0 \dyad \tena d_j^0 \big)\times \big[\big(\tena a^\alpha\!\cdot\!\tena d^0_\beta\big)\, \tena d^0_\beta\big] \\
& = \big(\tena a^\alpha\!\cdot\!\tena d^0_\beta\big) \big[      e_{ijs}\, K_{s\alpha}\,\tena d_i^0 \dyad \big(\tena d_j^0   \times  \tena d^0_\beta\big)\big] = \big(\tena a^\alpha\!\cdot\!\tena d^0_\beta\big) \big[    e_{ijs}\, e_{j\beta m}\, K_{s\alpha}\,\tena d_i^0 \dyad \tena d_m^0 \big] \\
& =  \big(\tena a^\alpha\!\cdot\!\tena d^0_\beta\big) \big[  \big( \delta_{im} \,\delta_{s\beta} - \delta_{i\beta} \,\delta_{sm}\big)  K_{s\alpha}\,\tena d_i^0 \dyad \tena d_m^0 \big] \\
& = \big(\tena a^\alpha\!\cdot\!\tena d^0_\beta\big) \big[  - K_{s\alpha}\,\tena d_\beta^0 \dyad \tena d_s^0 +  K_{\beta\alpha}\,\tena d_i^0 \dyad \tena d_i^0 \big]  \\
& = - K_{i\alpha} \big[ \big(\tena a^\alpha\!\cdot\!\tena d^0_\beta\big)\,\tena d_\beta^0 \big]\dyad \tena d_i^0 +  K_{\beta\alpha}\big(\tena d^0_\beta\!\cdot\!\tena a^\alpha\big)\,\tena 1_3 \big] \\
& = - \big(K_{i\alpha}  \, d_i^0\dyad \tena a^\alpha\big)^T + \mathrm{tr}\big(K_{i\alpha}  \, d_i^0\dyad \tena a^\alpha\big)\,\tena 1_3 = - \tenb K_e^T+ \big(\mathrm{tr}\,\tenb K_e\big)\tena 1_3\,,
\end{alignat*}
which shows that \eqref{e44}$_1$ holds true. Applying the trace operator to equation \eqref{e44}$_1$ we get $\mathrm{tr}\,\tenb K_e = \frac12\, \mathrm{tr}\,\tenb D_e\,$. Inserting this into \eqref{e44}$_1$ we obtain \eqref{e44}$_2\,$. The proof is complete.
\begin{remark}\label{rem9}
As a consequence of relations \eqref{e44} we deduce the relations between the norms, traces, symmetric and skew-symmetric parts of the two tensors   in the forms
\begin{alignat}{2}\label{e45}
  &  \|\tenb D_e \|^2 = \|\tenb K_e \|^2 + \big(\mathrm{tr}\,\tenb K_e\big)^2,\qquad \|\tenb K_e \|^2  = \|\tenb D_e \|^2 - \dfrac14\, \big(\mathrm{tr}\,\tenb D_e\big)^2,\\
  &  \mathrm{tr}\,\tenb D_e =  2\, \mathrm{tr}\,\tenb K_e\,,\quad \mathrm{skew}\,\tenb D_e =   \mathrm{skew}\,\tenb K_e\,,\quad \mathrm{dev_3sym}\,\tenb D_e =  - \mathrm{dev_3sym}\,\tenb K_e\,.\notag
\end{alignat}
Indeed the relations \eqref{e45} can be easily proved if we apply the operators $\mathrm{tr}$, $\|\cdot \|$, $\mathrm{skew}$, $\mathrm{dev}_3$, and $\mathrm{sym}$ to the equation \eqref{e44}$_1\,$.
In view of \eqref{e45}$_1\,$ and $\big( \mathrm{tr}\,\tenb K_e\big)^2\leq 3\,\|\tenb K_e \|^2$, we obtain the estimate
\begin{align} \label{e45,5}
\|\tenb K_e \|\leq \|\tenb D_e \|\leq 2\,\|\tenb K_e \|.
\end{align}
\end{remark}
In what follows, we write the components of the tensors $\tenb K_e$ and   $\tenb D_e\,$. To this aim, we use the relations
\begin{alignat}{2}\label{e46}
    \tenb Q_e^T\tenb Q_{e,\alpha} & =  \big( \tena d_i^0\dyad\tena d_i  \big)\big( \tena d_{k,\alpha}\dyad \tena d_k^0 + \tena d_{k}\dyad \tena d_{k,\alpha}^0  \big) \\
    & = \big( \tena d_i  \cdot \tena d_{k,\alpha}\big) \tena d_i^0\dyad\tena d_k^0 +  \tena d_i^0\dyad\tena d_{i,\alpha}^0 = \big( \tena d_i  \cdot \tena d_{k,\alpha} - \tena d^0_i  \cdot \tena d^0_{k,\alpha}\big) \tena d_i^0\dyad\tena d_k^0\,,\notag
\end{alignat}
which can be proved in the same way as equation \eqref{e14}. We compute the axial vector of the skew-symmetric tensor \eqref{e46} and find (similar to \eqref{e15})
\begin{align}\label{e47}
       \mathrm{axl}\big( \tenb Q_e^T\tenb Q_{e,\alpha}  \big) =  -\frac12 \, e_{ijk}\, \big(  \tena d_j\cdot\tena d_{k,\alpha} - \tena d^0_j\cdot\tena d^0_{k,\alpha} \big)\, \tena d^0_i\,.
\end{align}
By virtue of \eqref{e47} the definition \eqref{e41}  yields
\begin{alignat}{3}\label{e48}
      \tenb K_e\,  = & \, \frac12 \, e_{ijk}\, \big(  \tena d_{j,\alpha}\cdot\tena d_{k} - \tena d^0_{j,\alpha}\cdot\tena d^0_{k} \big)\, \tena d^0_i\dyad \tena a^\alpha\\
        = &  \,\big(  \tena d_{2,\alpha}\cdot\tena d_{3} - \tena d^0_{2,\alpha}\cdot\tena d^0_{3} \big)\, \tena d^0_1\dyad \tena a^\alpha + \big(  \tena d_{3,\alpha}\cdot\tena d_{1} - \tena d^0_{3,\alpha}\cdot\tena d^0_{1} \big)\, \tena d^0_2\dyad \tena a^\alpha \notag\\
        & + \big(  \tena d_{1,\alpha}\cdot\tena d_{2} - \tena d^0_{1,\alpha}\cdot\tena d^0_{2} \big)\, \tena d^0_3\dyad \tena a^\alpha ,\notag
\end{alignat}
which gives the components $K_{i\alpha}$ of the shell bending-curvature tensor $\tenb K_e$ in the tensor basis $\{ \tena d_i^0\dyad\tena a^\alpha \}$.

For the components of  $\tenb D_e\,$, we insert the relation \eqref{e46} in the equation \eqref{e43}
\begin{align*}
\tenb D_e = - \big( \tena d_i  \cdot \tena d_{k,\alpha} - \tena d^0_i  \cdot \tena d^0_{k,\alpha}\big)\big(  \tena d_i^0\dyad\tena d_k^0\big)\times \tena a^\alpha.
\end{align*}
Using that $\,\tena d_k^0 \times \tena a^\alpha= \tena d_k^0\times \big[\big(\tena a^\alpha\!\cdot\!\tena d^0_\beta\big)\, \tena d^0_\beta\big]= \big(\tena a^\alpha\!\cdot\!\tena d^0_\beta\big)\,e_{k\beta j}\,\tena d_j^0\,$, we obtain
\begin{align}\label{e49}
      \tenb D_e = e_{jk\beta}\big( \tena d_{i,\alpha}  \cdot \tena d_{k} - \tena d^0_{i,\alpha}  \cdot \tena d^0_{k}\big)\big(\tena a^\alpha\!\cdot\!\tena d^0_\beta\big)\, \tena d_i^0\dyad\tena d_j^0\,,
\end{align}
which shows the components   of the shell dislocation density tensor  in the tensor basis $\{ \tena d_i^0\dyad\tena d_j^0 \}$.

\section{Remarks and discussion}\label{sec:5}

Herein we present some other ways to express the shell dislocation density tensor, the shell bending-curvature tensor and discuss their close relationship.

\begin{remark}\label{rem10}
It is sometimes useful to express the components of the shell dislocation density tensor $\tenb D_e\,$ in the tensor basis $\{ \tena a^i\dyad\tena a_j \}$. If we multiply the relation  \eqref{e44}$_2$ with $\tena n_0$ and take into account that $\tenb K_e  \tena n_0 = \tena 0$, then we find $\tena 0 = - \tenb D^T_e \tena n_0 + \frac12\,\big( \mathrm{tr}\,\tenb D_e\big)\tena n_0 \,  $, which means
\begin{align*}
\tena n_0 \, \tenb D_e\,  =\, \frac12\,\big( \mathrm{tr}\,\tenb D_e\big)\tena n_0\,.
\end{align*}
It follows that the components of $\tenb D_e\,$ in the directions $\tena n_0\dyad\tena a_\alpha\,$ are zero, i.e. $\tenb D_e\,$ has the structure
\begin{align}\label{e50}
      \tenb D_e = \tenb D_\| + D_{\alpha 3}\, \tena a^\alpha\dyad\tena n_0 + \frac12\,\big( \mathrm{tr}\,\tenb D_e\big)\tena n_0\dyad\tena n_0\,,
\end{align}
where $\,\tenb D_\| = \tenb D_e\, \tena a = D^{\alpha\beta}\tena a_\alpha\dyad\tena a_\beta =  D_{\alpha\,\cdot}^{\,\,\,\,\beta}\,\tena a^\alpha\dyad\tena a_\beta $ is the planar part of $\tenb D_e\,$ (the part in the tangent plane).
If we insert \eqref{e50} into \eqref{e44}$_1$ and use $ \,\frac12\, \mathrm{tr}\,\tenb D_e = \mathrm{tr}\,\tenb K_e\,$, we get
\begin{align*}
\tenb D_\| + D_{\alpha 3}\, \tena a^\alpha\dyad\tena n_0 +  \big( \mathrm{tr}\,\tenb K_e\big)\tena n_0\dyad\tena n_0 = -K_{i\alpha}\,\tena a^\alpha\dyad\tena d_i^0 +  \big( \mathrm{tr}\,\tenb K_e\big) \big( \tena a +\tena n_0\dyad\tena n_0 \big),
\end{align*}
which implies (in view of \eqref{e48}) that
\begin{align}\label{e51}
     D_{\alpha 3} = -K_{3\alpha} = \tena d_{1}\cdot\tena d_{2,\alpha} - \tena d^0_{1}\cdot\tena d^0_{2,\alpha}\qquad\!\text{and}\qquad\! \tenb D_\| = - \big(\tenb K_\|\big)^T+ \big(\mathrm{tr}\,\tenb K_e\big)\tena a,
\end{align}
where  $\,\tenb K_\| \,= \tena a\,\tenb K_e\,  = K_{\beta\alpha}\tena d^0_\beta\dyad\tena a^\alpha $ is the planar part of $\tenb K_e\,$.
\end{remark}
\begin{remark}\label{rem11}
We observe that between the planar part $\,\tenb D_\| \,$ of $\,\tenb D_e\,$ and the planar part $\,\tenb K_\| \,$ of $\,\tenb K_e\,$ there exists a special relationship. The tensor $\,\tenb D_\| \,$ is the cofactor of the tensor $\,\tenb K_\| \,$. Let us explain this in more details: for any planar tensor $ \tenb S = S^\alpha_{\,\cdot\,\beta}\,\tena a_\alpha\dyad\tena a^\beta $ we introduce the transformation
\begin{align}\label{e52}
T( \tenb S ) =  - \tenb S^T +  \big( \mathrm{tr}\,\tenb S\big)\,\tena a\,.
\end{align}
One can prove that this transformation has the properties
\begin{align}\label{e53}
T\big(T( \tenb S)\big) = \tenb S\qquad\text{and}\qquad T( \tenb S ) =  - \tena c \,\tenb S\,   \tena c\,,
\end{align}
where the alternator $\tena c$ is defined in \eqref{e28}. Moreover, in view of \eqref{e53}$_2$ and \eqref{e28} we can write $ T( \tenb S )$ in the tensor basis $\{ \tena a^\alpha\dyad\tena a_\beta \}$  as follows
\begin{align*}
T( \tenb S ) = S^2_{\,\cdot\, 2}\,\tena a^1\dyad\tena a_1 - S^2_{\,\cdot\, 1}\,\tena a^1\dyad\tena a_2 - S^1_{\,\cdot\, 2}\,\tena a^2\dyad\tena a_1 + S^1_{\,\cdot\, 1}\,\tena a^2\dyad\tena a_2\,,
\end{align*}
which shows that the $2\times 2$ matrix of the components of $T(\tenb S)$ in the basis $\{ \tena a^\alpha\dyad\tena a_\beta \}$ is the cofactor of the matrix of components of $\tenb S$ in the basis $\{ \tena a_\alpha\dyad\tena a^\beta \}$, since
\begin{align*}
\left( \begin{array}{cc}
         S^2_{\,\cdot\, 2} & \,\,- S^2_{\,\cdot\, 1} \vspace{4pt}\\
         - S^1_{\,\cdot\, 2} & \,\,S^1_{\,\cdot\, 1}
       \end{array}
  \right) = \mathrm{Cof}
  \left( \begin{array}{cc}
         S^1_{\,\cdot\, 1} & \,\,S^1_{\,\cdot\, 2}  \vspace{4pt}\\
         S^2_{\,\cdot\, 1} & \,\, S^2_{\,\cdot\, 2}
       \end{array}
  \right).
\end{align*}
If the tensor $\tenb S$ is invertible, then from the Cayley-Hamilton relation  $\big( \tenb S^T\big)^2 - \big( \mathrm{tr}\,\tenb S\big)\tenb S^T + \mathrm{det} \tenb S = \tena 0$ and \eqref{e52} we deduce
\begin{align}\label{e54}
 T( \tenb S) =  - \tenb S^T +  \big( \mathrm{tr}\,\tenb S\big)\,\tena a = \big( \mathrm{det} \tenb S\big)\,\tenb S^{-T} =:  \mathrm{Cof}\big( \tenb S\big).
\end{align}
In our case, for the shell bending-curvature tensor $\tenb K_e$  we have $ \mathrm{tr}\, \tenb K_e = \mathrm{tr}\big(\tena a \tenb K_e\big) = \mathrm{tr}\big(  \tenb K_\|\big)  $, in view of \eqref{e48}. Then, the relation \eqref{e51}$_2$ yields
\begin{align*}
\tenb D_\| = - \big(\tenb K_\|\big)^T+ \big(\mathrm{tr}\,\tenb K_\|\big)\,\tena a\,.
\end{align*}
Using the relations \eqref{e52}-\eqref{e54} we see that $\,\tenb D_\| \,$ is the image of $\,\tenb K_\| \,$ under the transformation $T$, so that it holds
\begin{align}\label{e55}
\tenb D_\| = T\big(\tenb K_\|\big) =  - \tena c \,\big(\tenb K_\|\big)\,   \tena c = \mathrm{Cof}\big( \tenb K_\|\big),\\
\tenb K_\| = T\big(\tenb D_\|\big) =  - \tena c \,\big(\tenb D_\|\big)\,   \tena c = \mathrm{Cof}\big( \tenb D_\|\big). \notag
\end{align}
From  \eqref{e50}, \eqref{e51} we can write
\begin{align}\label{e56}
\tenb D_e = \mathrm{Cof}\big( \tenb K_\|\big) - K_{ 3\alpha}\, \tena a^\alpha\dyad\tena n_0 +  \big( \mathrm{tr}\,\tenb K_\|\big)\tena n_0\dyad\tena n_0\,,
\end{align}
which expresses once again the close relationship between the shell dislocation density tensor $\tenb D_e\,$ and the shell bending-curvature tensor $\tenb K_e\,$.
\end{remark}
\begin{remark}\label{rem12}
The shell bending-curvature tensor $\tenb K_e\,$ can also be expressed in terms of the directors $\tena d_i\,$. In this respect, an analogous relation to the formula \eqref{e24} for the wryness tensor (see Remark \ref{rem5}) holds
\begin{align}\label{e57}
 \tenb K_e \, = \, \,\frac12\, \big[  \tenb Q_e^T \big(\tena d_{i}\times  \tena d_{i,\alpha}\big) -    \tena d^0_{i}\times  \tena d^0_{i,\alpha}   \big]\dyad \tena a^\alpha .
\end{align}
To prove \eqref{e57}, we write the two terms in the brackets in the following form
\begin{alignat*}{2}
\tenb Q_e^T \big(\tena d_{i}\times  \tena d_{i,\alpha}\big) & = \big(\tena d_{k}^0\dyad  \tena d_{k}\big)\big(\tena d_{i}\times  \tena d_{i,\alpha}\big) =  \big[\tena d_{k}\cdot\big(\tena d_{i}\times  \tena d_{i,\alpha}\big)\big]\, \tena d_{k}^0 \\
& = \big[ \tena d_{i,\alpha} \cdot\big(\tena d_{k} \times \tena d_{i} \big)\big]\, \tena d_{k}^0 = e_{kij}\,\big(\tena d_{i,\alpha}\cdot\tena d_{j}\big)\, \tena d_{k}^0
\end{alignat*}
and similarly
\begin{align*}
\tena d^0_{i}\times  \tena d^0_{i,\alpha}  =   \big[\tena d_{k}^0\cdot\big(\tena d^0_{i}\times  \tena d^0_{i,\alpha}\big)\big]\, \tena d_{k}^0 = \big[ \tena d^0_{i,\alpha} \cdot\big(\tena d^0_{k} \times \tena d^0_{i} \big)\big]\, \tena d_{k}^0 = e_{kij}\,\big(\tena d^0_{i,\alpha}\cdot\tena d^0_{j}\big)\, \tena d_{k}^0.
\end{align*}
Inserting the last two relations into \eqref{e57} we obtain
\begin{align*}
\tenb K_e \, = \, \,\frac12\, e_{ijk}\,  \big[\big(  \tena d_{j,\alpha}\cdot\tena d_{k} \big)\, \tena d^0_i -  \big(\tena d^0_{j,\alpha}\cdot\tena d^0_{k} \big)\, \tena d^0_i  \big]\dyad \tena a^\alpha,
\end{align*}
which holds true, by virtue of \eqref{e48}. Thus, \eqref{e57} is proved.

We can put the relation \eqref{e57} in the form
\begin{align}
 & \tenb K_e \, = \, \tenb Q_e^T \,\tenb\omega\qquad\text{where we define}   \label{e58}\\
 & \tenb\omega = \tenb\omega_\alpha\dyad \tena a^\alpha\qquad\text{with}\qquad  \tenb\omega_\alpha = \,\frac12\, \big[ \tena d_{i}\times  \tena d_{i,\alpha} -  \tenb Q_e \big( \tena d^0_{i}\times  \tena d^0_{i,\alpha}\big) \big]. \label{e59}
\end{align}
If we compare the relations \eqref{e58}  and the definition \eqref{e41}, we derive
\begin{align*}
\tenb\omega_\alpha = \tenb Q_e\,\mathrm{axl}\big( \tenb Q_e^T\tenb Q_{e,\alpha}  \big) = \mathrm{axl}\big( \tenb Q_{e,\alpha} \,\tenb Q_e^T  \big).
\end{align*}
Then, from \eqref{e13} we deduce $\,\,\,
 \tenb Q_{e,\alpha} \,\tenb Q_e^T = \tenb\omega_\alpha \times \tena 1_3
\,\,\,$
and by multiplication with $\,\tenb Q_e\,$ we find
\begin{align}\label{e60}
\tenb Q_{e,\alpha} \,= \, \tenb\omega_\alpha \times \tenb Q_e\,\,,\qquad \alpha=1,2.
\end{align}
Thus, the equations \eqref{e58}, \eqref{e59} can be employed for an alternative definition of the shell bending-curvature tensor, namely
\begin{align}\label{e61}
  \tenb K_e \, = \, \tenb Q_e^T \,\tenb\omega\,,\qquad\text{where}\qquad \! \tenb\omega = \tenb\omega_\alpha\dyad \tena a^\alpha \!\!\qquad\text{and}\qquad\!\! \tenb Q_{e,\alpha} \,= \, \tenb\omega_\alpha \times \tenb Q_e\,.
\end{align}
This is the counterpart of the relations \eqref{e20}, \eqref{e21} for the wryness tensor in the three-dimensional theory of Cosserat continua. The relations \eqref{e61} were used to define the corresponding shell bending-curvature tensor,  e.g., in \cite{Altenbach04,Zhilin06}.
\end{remark}

\begin{remark}\label{rem13}
As shown by relations \eqref{e0,7} for the three-dimensional case, one can introduce the elastically stored shell energy density $W$ as a function of the  shell strain tensor and the shell dislocation density tensor
\begin{align}\label{e62}
W=W\big(\tenb E_e \,,\, \tenb D_e \big).
\end{align}
If \eqref{e62} is assumed to be a quadratic convex and coercive function, then the existence of solutions to the minimization problem of the total energy functional for Cosserat shells can be proved in a similar manner as in Theorem 6 of \cite{Birsan-Neff-MMS-2014}. In the proof, one should employ decisively the estimate \eqref{e45,5} and the expressions of the shell dislocation density tensor $\, \tenb D_e\,$ established in the previous sections.
\end{remark}


\end{document}